\documentclass[letterpaper, 10 pt, conference]{ieeeconf}  
\IEEEoverridecommandlockouts                              
\usepackage{epsfig} 
\usepackage{graphicx}
\usepackage{amsmath}
\usepackage{array}
\usepackage{amssymb}
\usepackage{verbatim}
\usepackage{cleveref}
\usepackage{color}
\usepackage{relsize}
\usepackage{tikz}
\usepackage{fancyhdr}
\usepackage{graphics}
\usepackage{setspace}
\usepackage{epstopdf}
\usepackage{multirow}
\usepackage{lipsum}
\usepackage{textcomp}
\usepackage{mathtools}
\usepackage{cite}
\usepackage{setspace}
\usepackage{booktabs}
\usepackage{breqn}
\usepackage{algorithm}
\usepackage[noend]{algpseudocode}
\usepackage{flushend}
\usepackage{multicol}
\usepackage{float}
\usepackage{mathtools,bm}
\usepackage[noend]{algpseudocode}
\usepackage{booktabs}
\usepackage{dblfloatfix}
\usepackage{mathrsfs}
\usepackage{mathtools}
\usepackage{leftidx}
\usepackage[utf8]{inputenc}
\usepackage[english]{babel}

\fancypagestyle{firstpage}{\lhead{\vspace{-0.75 cm} \small \textbf{{\fontfamily{cmss}\selectfont
58th Conference on Decision and Control (CDC)\\December 11--13, 2019, Nice, France}}}}

\title{\large \bf
Robust Hierarchical MPC for Handling Long Horizon Demand Forecast Uncertainty\\with Application to Automotive Thermal Management$^*$}

\author{Mohammad Reza Amini$^{1}$, Ilya Kolmanovsky$^{2}$, and Jing Sun$^{1}$
\thanks{*This paper is based upon the work supported by the United States Department of Energy (DOE) under award No. DE-AR0000797.}
\thanks{$^{1}$M.R. Amini and J. Sun are with the Department of Naval Architecture \& Marine Engineering, University of Michigan, Ann Arbor, MI 48109 USA. 
Emails: {\tt\small \{mamini,jingsun\}@umich.edu}}%
\thanks{$^{2}$I. Kolmanovsky is with  the Dept. of Aerospace Engineering, University of Michigan, Ann Arbor, MI 48109 USA. Email: {\tt\small ilya@umich.edu}}%
}

\renewcommand{\baselinestretch}{1}

\begin{document}

\maketitle
\thispagestyle{firstpage}

\begin{abstract}
This paper presents a robust hierarchical MPC (H-MPC) for dynamic systems with slow states subject to demand forecast uncertainty. The H-MPC has two layers: (i) the scheduling MPC at the upper layer with a relatively long prediction/planning horizon and slow update rate, and (ii) the piloting MPC at the lower layer over a shorter prediction horizon with a faster update rate. The scheduling layer MPC calculates the optimal slow states, which will be tracked by the piloting MPC, while enforcing the system constraints according to a long-range and approximate prediction of the future demand/load, e.g., traction power demand for driving a vehicle. In this paper, to enhance the H-MPC robustness against the long-term demand forecast uncertainty, we propose to use the high-quality preview information enabled by the connectivity technology over the short horizon to modify the planned trajectories via a constraint tightening approach at the scheduling layer. Simulation results are presented for a simplified vehicle model to confirm the effectiveness of the proposed robust H-MPC framework in handling demand forecast uncertainty.
 \end{abstract}

\section{INTRODUCTION}
For dynamic systems with multiple time scales, centralized optimization could lead to non-uniform performance as a centralized controller inherently assumes that all the system outputs are equally reactive~\cite{medagoda2014multiple}. However, the outputs of such dynamic systems tend to react and respond at different rates. This is a common observation for dynamic systems with slow and storage states, e.g., microgrids~\cite{clarke2018hierarchical}, vehicle~\cite{AminiACC19,AminiCDC18,Amini_TCST_2019,koeln2018robust} and aircraft~\cite{koeln2019hierarchical} thermal management systems.

One approach to address the aforementioned issue in a centralized scheme is to select a relatively long prediction horizon to ensure satisfactory performance of the slow states.~While this approach could fulfill the overall performance requirements, it has several significant issues for practical applications:
\begin{itemize}
    \item Firstly, it will significantly increase the computational complexity. For fast states, the computational complexity is often not justifiable as the information downstream of the horizon has minimal influence on the optimal evolution of the fast state trajectories.
    \item Secondly, for most systems operating in a dynamic environment, accurate prediction of the demand/load profile over a long prediction horizon is not feasible. For example, vehicle speed cannot be predicted accurately over a long prediction horizon based on vehicle-to-vehicle (V2V) and vehicle-to-infrastructure (V2I) information~\cite{Amini_TCST_2019,AminiCDC18}. 
\end{itemize}  
The common approach exploited in the literature~\cite{koeln2018robust,barcelliy2010hierarchical,scattolini2007hierarchical} for systems with multi-timescales is the hierarchical control with multiple layers. The hierarchical framework allows for optimization over different planning horizons at each layer based on the requirements of different system outputs. It also helps to significantly reduce the computational cost compared to its centralized counterpart~\cite{AminiCDC18,koeln2018robust,scattolini2009architectures}. 

Whether centralized or hierarchical, the uncertainty associated with the long-term prediction of the demand/load is the common problem for systems with slow dynamics. As the optimization often pushes the system to work at its limit for best efficiency, it is important to be able to enforce the constraints robustly in the presence of demand forecast uncertainty. This is because for slow states, once the violation of the constraint occurs, the controller has to put more efforts to compensate for this violation, thereby reducing the overall system efficiency. Unlike the most common techniques in the literature on the robust design of hierarchical and distributed MPC~\cite{koeln2018robust,scattolini2007hierarchical,picasso2010mpc} which mainly focus on robustifying the controller under the presence of external unknown disturbances, the robustness to long horizon demand forecast uncertainty has not been fully studied. 

From a new perspective, in this paper, we consider the long horizon demand forecast uncertainty challenge which plays a critical role in planning and optimization of slow/storage states with wide application in vehicle and aircraft thermal management, building HVAC control, and power grid optimization. To address this challenge, we exploit the H-MPC optimization framework as the hierarchical architecture provides a unique capability of using the demand forecast with different accuracies over different prediction horizons. Our proposed H-MPC has two layers as scheduling and piloting layers. The scheduling layer at the top computes the optimal control input at a slow update rate over a relatively long prediction horizon. The output of the scheduling layer is the optimal trajectories of the ``slow'' states of interest according to an approximate prediction of the demand over the long planning horizon. In the lower layer, piloting layer, over a much shorter prediction horizon and using a faster update rate, an MPC is designed to follow the scheduled slow state trajectories coming from the upper layer while fulfilling the overall system performance requirements.

It is envisioned that over a relatively short prediction horizon, the demand/load profile can be accurately predicted. As an example for connected and automated vehicles, the V2V/V2I communications can be used to predict the future vehicle speed and the associated traction power demand over a short time window accurately~\cite{Zhen2018Traj}. We propose to use this short horizon accurate demand preview information to predict the system state evolution over the prediction horizon of the piloting layer MPC while tracking the planned trajectory coming from the scheduling layer MPC. Computing the evolution of the system state with a faster update rate when the accurate demand preview is accessible allows for detecting the constraint violation proactively over the short prediction horizon. Particularly, we can predict the deviation of the slow state from the planned trajectory and realize the possibility of constraint violation over the short prediction horizon. It is noted that, assuming a perfect model, the deviation in the reference tracking and constraint violation at the piloting layer are mainly attributed to the mismatch between the actual demand and the approximate preview incorporated for planning at the scheduling layer. Once these deviations and possible constraint violation are predicted over the shorter prediction horizon, the constraint set at the scheduling layer is re-computed and tightened to make the H-MPC prepared for the upcoming situation during which the chances of constraint violation are high. 

The main contribution of this paper is to develop a robust H-MPC for a class of dynamic systems with slow states subject to demand preview uncertainty. We propose to utilize short and long-range demand preview information to compute robust constraint set to enhance the conventional H-MPC robustness against the demand preview uncertainty while improving the overall system efficiency. The performance of the proposed robust H-MPC is studied through simulations for a vehicle thermal management case study.

\vspace{-0.1cm} 
\section{Hierarchical MPC Formulation}\label{sec:Sec2}\vspace{-0.05cm} 
We consider the following discrete-time linear time-invariant system,
\vspace{-0.25cm}
\begin{gather}\label{eq:ACC18_Eq1}
\boldsymbol{x}_{k+1}={A}\boldsymbol{x}_{k}+{B}_1\tilde{\boldsymbol{u}}_{k}+{B}_2\hat{\boldsymbol{u}}_{k}, 
\end{gather}
where $\boldsymbol{x}=[x_1,\cdots,x_n]^\intercal \in \mathbb{R}^n$, $\tilde{\boldsymbol{u}} \in \mathbb{R}^m$, and $\hat{\boldsymbol{u}} \in \mathbb{R}^h$ are the system state, control input, and demand/load on the system, respectively, and ${A}\in\mathbb{R}^{n\times n}$, ${B}_1\in\mathbb{R}^{n\times m}$, and ${B}_2\in\mathbb{R}^{n\times h}$. The integer $k$ denotes the time step and the sampling period of (\ref{eq:ACC18_Eq1}) is assumed to be $T$.

We further assume that different states of the system ($x_{i_x},~i_x=1,\cdots,n)$ respond at different rates. Depending on the relative time scales of different states, they can be categorized into ``fast'' and ``slow'' states: $\boldsymbol{x}=[({\boldsymbol{x}^{\text{\textbf{fast}}}})^{\intercal},{(\boldsymbol{x}^{\text{\textbf{slow}}}})^{\intercal}]^{\intercal}$.~Note that $\tilde{\boldsymbol{u}}$ contains the physical control inputs to the system that are optimized to regulate the overall performance. On the other hand, $\hat{\boldsymbol{u}}$ is the vector of non-adjustable inputs to the system representing measured external loads. Since $\hat{\boldsymbol{u}}$ is a measured input to the system, it affects the system states, and the knowledge about $\hat{\boldsymbol{u}}$ can be incorporated as a preview for optimization of $\tilde{\boldsymbol{u}}$. In automotive applications, traction power demand ($P_{trac}$) for driving a vehicle is an example of $\hat{\boldsymbol{u}}$.

The state and control input of the linear system in (\ref{eq:ACC18_Eq1}) are constrained: $\boldsymbol{x} \in\mathcal{{X}}$, $\tilde{\boldsymbol{u}} \in\mathcal{{U}}$, where~ 
$\mathcal{{X}}=\{\boldsymbol{x}_k|P_x\boldsymbol{x}_k\leq \boldsymbol{q}_x\}$ and $\mathcal{{U}}=\{\tilde{\boldsymbol{u}}_k|P_u\tilde{\boldsymbol{u}}_k\leq \boldsymbol{q}_u\}$ are polyhedral sets.%

\vspace{-0.1cm}
\subsection{Single-Layer MPC (S-MPC) Formulation}
For the system (\ref{eq:ACC18_Eq1}), we first consider a centralized single-layer MPC (S-MPC) defined according to the following finite-horizon optimization problem:\vspace{-0.15cm} %
\begin{gather} \label{eqn:S_MPC_gen}
\begin{aligned}
& \underset{{\tilde{\boldsymbol{U}}}_k}{\text{min}} & & \sum_{j=0}^{N} \ell_{\text{S-MPC}}(k+j|k), \\
& \text{s.t.}
& & \boldsymbol{x}(k+j+1|k)=A\boldsymbol{x}(k+j|k)+B_1\tilde{\boldsymbol{u}}(k+j|k)\\
& 
& & ~~~~~~~~~~~~~~~~~~+B_2\hat{\boldsymbol{u}}(k+j|k),\\
& 
& & \boldsymbol{x}(k+j|k)\in \mathcal{X},~\tilde{\boldsymbol{u}}(k+j|k)\in \mathcal{U},\\
& 
& & \boldsymbol{x}(k|k)=\boldsymbol{x}_k,
\end{aligned}
\end{gather}
%
where $\ell_{\text{S-MPC}}$ is the S-MPC stage cost, $\tilde{\boldsymbol{U}}_k=[\tilde{\boldsymbol{u}}(k|k)^{\intercal},\cdots,\tilde{\boldsymbol{u}}(k+N-1|k)^{\intercal}]^\intercal$ is the control input sequence, and $\hat{\boldsymbol{U}}_k=[\hat{\boldsymbol{u}}(k|k)^{\intercal},\cdots,\hat{\boldsymbol{u}}(k+N-1|k)^{\intercal}]^\intercal$ is the sequence of the ``predicted'' demand preview. The cost function $\ell_{\text{S-MPC}}$ is defined as~\cite{koeln2018two}:\vspace{-0.15cm}
\begin{gather}\label{eq:SMPC_cost_gen}
    \ell_{\text{S-MPC}}=\left\Vert \boldsymbol{r}({k+j|k})-\boldsymbol{r}^{d}({k+j}) \right\Vert^2_\Lambda
\end{gather}
where $\boldsymbol{r}_k=E\boldsymbol{x}_k+F\tilde{\boldsymbol{u}}_k$, and $\boldsymbol{r}^{d}\in \mathbb{R}^{n_r}$ represents the desired operation of the system. When using the solution of the optimization problem (\ref{eqn:S_MPC_gen}) to control the system (\ref{eq:ACC18_Eq1}), the MPC feedback law defined by $\tilde{\boldsymbol{u}}(k|k)$ inherently assumes that all the system outputs are equally reactive~\cite{medagoda2014multiple}, however, in the case considered here, the outputs of the system respond over different time scales.

\vspace{-0.1cm}
\subsection{Hierarchical MPC (H-MPC) Formulation}\vspace{-0.1cm}
The H-MPC considered in this paper has two layers~\cite{scattolini2007hierarchical}:
\begin{itemize}
    \item \textbf{scheduling layer} at the top with an MPC that computes the optimal control inputs and the ``slow'' state trajectories at a slower update rate of $T_s=\nu T$, where $\nu$ is a positive integer, over a relatively long prediction horizon ($H_s$). An approximate prediction of the demand sequence is incorporated at this layer.
    \item \textbf{piloting layer} with an MPC implemented over a much shorter prediction horizon ($H_p$) and using a faster update rate ($T$) to track the planned trajectories of the slow states from the scheduling layer. The piloting layer MPC has access to a more accurate prediction of the demand profile over $H_p$. 
\end{itemize} 

The scheduling layer MPC is based on the following finite-horizon optimization problem:\vspace{-0.2cm}
\begin{gather} \label{eqn:H_MPC_UL_gen}
\begin{aligned}
& \underset{\tilde{\boldsymbol{U}}^s_{k_s}}{\text{min}} & & \sum_{j=0}^{H_s} \ell_{\text{scheduling}}(k_s+j|k_s), \\
& \text{s.t.}
& & \boldsymbol{x}^s(k_s+j+1|k_s)=A^s\boldsymbol{x}^s(k_s+j|k_s)\\
& 
& & +B^s_1\tilde{\boldsymbol{u}}_s(k_s+j|k_s)+B^s_2\hat{\boldsymbol{u}}_s(k_s+j|k_s),\\
& 
& & \boldsymbol{x}^s(k_s+j|k_s)\in \mathcal{X},~\tilde{\boldsymbol{u}}^s(k_s+j|k_s)\in \mathcal{U},\\
& 
& & \boldsymbol{x}^s(k_s|k_s)=\boldsymbol{x}^s_{k_s},
\end{aligned}
\end{gather}
where $\boldsymbol{x}^s$ and $\tilde{\boldsymbol{U}}^s_{k_s}=[\tilde{\boldsymbol{u}}^s(k_s|k_s)^{\intercal},\cdots,\tilde{\boldsymbol{u}}^s(k_s+H_s-1|k_s)^{\intercal}]^\intercal$ are the state vector and control input sequence updated at the slower sampling time of $T_s=\nu T$, which is indexed by $k_s$. The sequence of the approximate demand profile $\hat{\boldsymbol{U}}^s_{k_s}=[\hat{\boldsymbol{u}}^s(k_s|k_s)^{\intercal},\cdots,\hat{\boldsymbol{u}}^s(k_s+H_s-1|k_s)^{\intercal}]^\intercal$ is also incorporated at the scheduling layer as the preview. The stage cost function ($\ell_{\text{scheduling}}$) of (\ref{eqn:H_MPC_UL_gen}) is the same as the one considered in (\ref{eq:SMPC_cost_gen}) while it is updated every $T_s$ second. The prediction model in (\ref{eqn:H_MPC_UL_gen}) is based on the down-sampled version of (\ref{eq:ACC18_Eq1}) with slower sampling rate of $T_s$ as follows:\vspace{-0.2cm}
\begin{gather}
    A^s=(A)^{\nu},~~B^s_1=\sum_{j=0}^{\nu-1}(A)^jB_1,~~B^s_2=\sum_{j=0}^{\nu-1}(A)^jB_2,
\end{gather}

Once the scheduling layer optimization problem is solved, based on the computed control sequence $\boldsymbol{U}^s_{k_s}$, the sequence of the scheduled slow states at the subsequent time steps over the horizon of the piloting layer MPC ($H_p$) can be calculated: $\leftidx{^*}{\boldsymbol{x}}^{\text{\textbf{slow}}}(k_s+1|k_s),\cdots,\leftidx{^*}{\boldsymbol{x}}^{\text{\textbf{slow}}}(k_s+H_p|k_s)$, where $\leftidx{^*}{\boldsymbol{x}}^{\text{\textbf{slow}}}$ is the vector of scheduled slow states of interest. Note that unlike the approach presented in~\cite{scattolini2007hierarchical,koeln2018two}, the horizon of the second layer, piloting layer ($H_p$), can be longer than $H_p>\nu$. In the next step, the sequence of $\leftidx{^*}{\boldsymbol{x}}^{\text{\textbf{slow}}}$ is passed to the piloting layer as a piecewise constant function to account for the faster sampling rate of the lower layer controller ($T$):\vspace{-0.15cm}
\begin{gather}
    \leftidx{^*}{\boldsymbol{X}}^{\text{\textbf{slow}}}(k|k_s)=[\leftidx{^*}{\boldsymbol{x}}^{\text{\textbf{slow}}}(k_s|k_s)^{\intercal},\leftidx{^*}{\boldsymbol{x}}^{\text{\textbf{slow}}}(k+1|k_s)^{\intercal},\cdots, \nonumber\\
    ~~~~~~~~~~\leftidx{^*}{\boldsymbol{x}}^{\text{\textbf{slow}}}(k+\nu-1|k_s)^{\intercal},\leftidx{^*}{\boldsymbol{x}}^{\text{\textbf{slow}}}(k_s+1|k_s)^{\intercal},\cdots,\nonumber\\
    ~~\leftidx{^*}{\boldsymbol{x}}^{\text{\textbf{slow}}}(k+2\nu-1|k_s)^{\intercal},\cdots,\leftidx{^*}{\boldsymbol{x}}^{\text{\textbf{slow}}}(k+H_p-1|k_s)^{\intercal}]^\intercal\label{eq:H_MPC_LL_Ref}
\end{gather}
where $\leftidx{^*}{\boldsymbol{X}}^{\text{\textbf{slow}}}$ is the sequence of the scheduled slow states to be tracked at the piloting layer. Note that the sequence of $\leftidx{^*}{\boldsymbol{X}}^{\text{\textbf{slow}}}$ between the slow updates are assumed to be constant, e.g., $\leftidx{^*}{\boldsymbol{x}}^{\text{\textbf{slow}}}(k_s|k_s)=\leftidx{^*}{\boldsymbol{x}}^{\text{\textbf{slow}}}(k+1|k_s)=\cdots=\leftidx{^*}{\boldsymbol{x}}^{\text{\textbf{slow}}}(k+\nu-1|k_s)$.
~The piloting layer MPC is based on the solution of the following optimal control problem:\vspace{-0.2cm}
\begin{gather} \label{eqn:H_MPC_LL_gen}
\begin{aligned}
& \underset{\tilde{\boldsymbol{U}}_k}{\text{min}} & & \sum_{j=0}^{H_p} \ell_{\text{piloting}}(k+j|k), \\
& \text{s.t.}
& & \boldsymbol{x}(k+j+1|k)=A\boldsymbol{x}(k+j|k)+B_1\tilde{\boldsymbol{u}}(k+j|k)\\
& 
& & ~~~~~~~~~~~~~~~~~~+B_2\Hat{\boldsymbol{u}}(k+j|k),\\
& 
& & \boldsymbol{x}(k+j|k)\in \mathcal{X},~\tilde{\boldsymbol{u}}(k+j|k)\in \mathcal{U},\\
& 
& & \boldsymbol{x}(k|k)=\boldsymbol{x}_{k}.
\end{aligned}
\end{gather}
%
where the stage cost $\ell_{\text{piloting}}$ is a combination of the original cost in (\ref{eq:SMPC_cost_gen}) and extra tracking terms to enforce the slow states ($\boldsymbol{x}^{\text{\textbf{slow}}}$) to track the scheduled slow state trajectories ($\leftidx{^*}{\boldsymbol{x}}^{\text{\textbf{slow}}}$) presented in (\ref{eq:H_MPC_LL_Ref}). ${\tilde{\boldsymbol{U}}}_k=[\tilde{\boldsymbol{u}}(k|k)^{\intercal},\cdots,\tilde{\boldsymbol{u}}(k+H_p-1|k)^{\intercal}]^\intercal$ is the control input sequence of the piloting layer MPC. It is assumed that the incorporated demand preview at the piloting layer (${\hat{\boldsymbol{U}}}_k=[\hat{\boldsymbol{u}}(k|k)^{\intercal},\cdots,\hat{\boldsymbol{u}}(k+H_p-1|k)^{\intercal}]^\intercal$) is more accurate than the demand preview used at the scheduling layer (${\hat{\boldsymbol{U}}}^s_{k_s}$). When implementing the H-MPC, at the end of each optimization iteration, $\tilde{\boldsymbol{u}}(k|k)$ is commanded to the system and the horizon is shifted by one time step ($T$).

\vspace{-0.15cm}
\section{Robust Hierarchical MPC}\label{sec:sec3}
%

\begin{figure}[h!]
	\begin{center}
		\includegraphics[width=1.0\columnwidth]{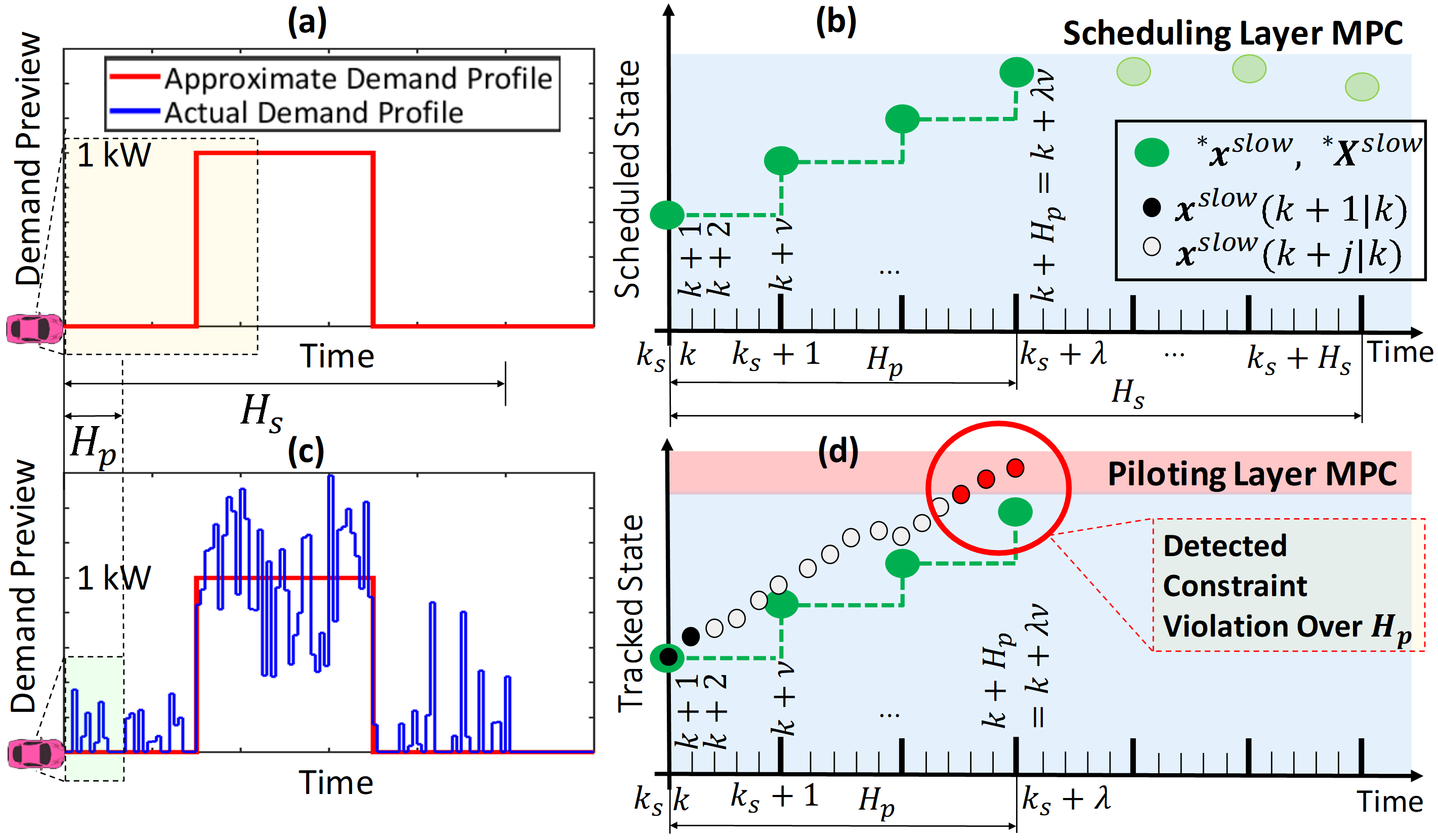} \vspace{-0.65cm}   
		\caption{Schematic of the constraint violation prediction to compute $\mathcal{X}^{\text{\textbf{robust}}}$. At time step $k$, the scheduling layer MPC plans the slow state trajectories (b) based on an approximate knowledge of the demand preview (a). The piloting layer MPC is designed to not only regulate the fast states but track the planned slow state trajectory while accessing the accurate demand preview at a faster update rate. The evolution of the slow states over the short prediction horizon is computed (d). If any deviation or constraint violation is predicted over $H_p$, the scheduling layer constraint set is tightened proactively.}\vspace{-0.55cm} 
		\label{fig:ErrorPredictionConcept} 
	\end{center}
\end{figure}

An intuitive approach for improving the constraint satisfaction at the piloting layer is to tighten the constraints at the scheduling layer. If the constraint tightening is done conservatively, the overall efficiency of the system may be degraded. To avoid unnecessary conservatism, in this paper, we propose to compute a robust constraint set (denoted by $\mathcal{X}^{\text{\textbf{robust}}}$) at each time step $k$. With a faster update rate and more accurate incorporated demand preview, the evolution of the slow state trajectories over the short horizon of the piloting layer MPC can be predicted while tracking the planned trajectories from the scheduling layer. This prediction of the slow state evolution is then used to compute $\mathcal{X}^{\text{\textbf{robust}}}$ via estimating the deviation in tracking the planned slow state trajectory and potential violation of the constraint over $H_p$. This concept is shown in Fig.~\ref{fig:ErrorPredictionConcept}. 

According to the accurate demand preview ($\hat{\boldsymbol{U}}_k=[\hat{\boldsymbol{u}}(k|k)^{\intercal},\cdots,\hat{\boldsymbol{u}}(k+H_p-1|k)^{\intercal}]^\intercal$), the evolution of the system states $\boldsymbol{x}$ over the piloting layer horizon $H_p$ can be predicted based on (\ref{eq:ACC18_Eq1}):
\vspace{-0.15cm}
\begin{gather}
    \boldsymbol{x}(k+j+1|k)=A\boldsymbol{x}(k+j|k)~~~~~~~~~~~~~~~~~~~~~~~~~~~~~\nonumber\\
    ~~~~~~~~~~~+B_1\tilde{\boldsymbol{u}}^{*}(k+j|k)+B_2\hat{\boldsymbol{u}}(k+j|k)\label{eq:X_evolution}
\end{gather}
where $j=0,\cdots,H_p-1$ and $\tilde{\boldsymbol{u}}^{*}(k+j|k)$ is the vector of the control inputs at $k+j$ computed based on the piloting layer optimal control problem (\ref{eqn:H_MPC_LL_gen}) solution at time step $k$. Based on the state predictions (\ref{eq:X_evolution}), the constraint violation of the slow states of interest can also be predicted at subsequent time steps, which will be used to tighten the state constraint set for the optimization problem at the scheduling layer. We denote the tightened constraint set over the piloting layer horizon at $k+j+1$, $j=0,\cdots,H_p-1$, by $\mathcal{X}_j$, which is computed as follows: 
\vspace{-0.15cm}
\begin{gather}\label{eq:X_robust_1}
   \mathcal{X}_{j}=\{\boldsymbol{x}(k|k)|~~~~~~~~~~~~~~~~~~~~~~~~~~~~~~~~~~~~~~~~~~~~~\\
   ~~~~P_x\boldsymbol{x}(k+j+1|k)\leq\big(\boldsymbol{q}_x-\bar{\mathbf{q}}_x(k+j+1|k)\big)\},\nonumber
\end{gather}
where $\boldsymbol{q}_x=[q_x^1,\cdots,q_x^{n_q}]^{\intercal}\in\mathbb{R}^{n_q}$ represents the state constraints, and is used in the definition of the original state constraint set $\mathcal{X}=\{\boldsymbol{x}_k|P_x\boldsymbol{x}_k\leq \boldsymbol{q}_x\}$. Note that if we assume $\boldsymbol{q}_x$ represents the upper and lower bounds on each state and each row of $P_x$ corresponds to one constraint, then non-zero elements of $P_x$ are either $1$ (upper limit) or $-1$ (lower limit). We then define $\bar{\mathbf{q}}_x(k+j+1|k)$ in (\ref{eq:X_robust_1}) as the vector of predicted constraint violation at $k+j+1$ based on the state prediction made at $k$ (\ref{eq:X_evolution}). If $q_x^{i_q}$ ($i_q=1,\cdots,n_q$) corresponds to $x_{i_x}$ ($i_x=1,\cdots,n$), then the $i_q^{\text{th}}$ element of $\bar{\mathbf{q}}_x$ at $k+j+1$, i.e., $\bar{{q}}_x^{i_q}(k+j+1|k)$, is calculated as follows:
\vspace{-0.15cm}
\begin{gather}\label{eq:X_robust_0}
   \bar{q}_x^{i_q}(k+j+1|k)=~~~~~~~~~~~~~~~~~~~~~~~~~~~~~~~~~~~~~~\\\begin{cases}0~~~~\text{\textbf{if}}~x_{i_x}~\text{is \textbf{fast}}\\0~~~~\text{\textbf{if}}~x_{i_x}(k+j+1|k)<q_x^{i_q}~\&~x_{i_x}~\text{is \textbf{slow}}\\\text{}~\vspace{-0.35cm}\\
   x_{i_x}(k+j+1|k)-q_x^{i_q}\\
   ~~~~~~\text{\textbf{if}}~x_{i_x}(k+j+1|k)\geq q_x^{i_q}~\&~x_{i_x}~\text{is \textbf{slow}}
   \end{cases}\nonumber
\end{gather}
Note that since we are interested in the slow states constraint enforcement, the elements of $\bar{\mathbf{q}}_x$ corresponding to fast states are set to zero in (\ref{eq:X_robust_0}).

Eventually and based on the computed $\mathcal{X}_{j}$s in (\ref{eq:X_robust_1}), the robust constraint set ($\mathcal{X}^{\text{\textbf{robust}}}$) to be implemented at the next time step $k+1$ at the scheduling layer is computed as:\vspace{-0.1cm}
\begin{gather}\label{eq:X_robust_2}
   \mathcal{X}^{\text{\textbf{robust}}}=\mathcal{X}\cap\mathcal{X}_{1}\cap\cdots \cap\mathcal{X}_{k+H_p-1},
\end{gather}
based on which, the scheduling layer MPC of the robust H-MPC is re-formulated as follows: \vspace{-0.15cm}
\begin{gather} \label{eqn:H_MPC_UL_gen_robust}
\begin{aligned}
& \underset{\tilde{\boldsymbol{U}}^s_{k_s}}{\text{min}} & & \sum_{j=0}^{H_s} \ell_{\text{scheduling}}(k_s+j|k_s), \\
& \text{s.t.}
& & \boldsymbol{x}^s(k_s+j+1|k_s)=A^s\boldsymbol{x}^s(k_s+j|k_s)\\
& 
& & +B^s_1\tilde{\boldsymbol{u}}_s(k_s+j|k_s)+B^s_2\hat{\boldsymbol{u}}_s(k_s+j|k_s),\\
& 
& & \boldsymbol{x}^s(k_s+j|k_s)\in \mathcal{X}^{\text{\textbf{robust}}},\\
& 
& & \tilde{\boldsymbol{u}}^s(k_s+j|k_s)\in \mathcal{U},~\boldsymbol{x}^s(k_s|k_s)=\boldsymbol{x}^s_{k_s},
\end{aligned}
\end{gather}
%

\hspace{-0.2cm}\textbf{{Remark}}: Note that one can propose to tighten the constraint whenever a violation of the constraint is detected at the piloting layer at time step $k$. This event-triggered ``passive'' constraint tightening approach can be incorporated in the H-MPC at the scheduling layer (\ref{eqn:H_MPC_UL_gen}) by replacing the original constraint set $\mathcal{X}$ with $\mathcal{X}^{\text{\textbf{passive}}}$, which is computed as follows: \vspace{-0.35cm}  
\begin{gather}\label{eq:X_passive}
   \mathcal{X}^{\text{\textbf{passive}}}=\{\boldsymbol{x}(k|k)|P_x\boldsymbol{x}(k|k)\leq (q_x-\bar{\mathbf{q}}_x(k|k))\},
\end{gather}
where $\bar{\mathbf{q}}_x(k|k)$ is the vector of measured (detected) constraint violation at each time step $k$.~
While $\mathcal{X}^{\text{\textbf{passive}}}$ is easy to compute, it leads to marginal improvement in the H-MPC robustness, see~\cite{AminiCDC18} for an example. This can be explained with respect to the slow dynamic characteristics which call for a long prediction horizon with enough lead time to ensure an effective constraint enforcement. In other words, if the H-MPC constraint set is tightened based on (\ref{eq:X_passive}) after the occurrence of the violation, compensating for the occurred violation in the subsequent time steps may not be effective as it takes extra effort from the controller, potentially reducing the system efficiency.

\section{Case Study: Vehicle Thermal Management}
\vspace{-0.1cm}

A simplified model of a vehicle is considered in this section to demonstrate the effectiveness of the proposed robust H-MPC. The main structure of the model is adopted from~\cite{koeln2018two} and a new state has been added to the model to represent the slow state. The model has four states, including the position of the vehicle ($x_1$), vehicle speed ($x_2$), on-board stored energy ($x_3$), e.g., battery state of charge, and thermal index of the vehicle energy storage ($x_4$), e.g., battery temperature. Compared to the first three states, the fourth state has a relatively slower dynamics. The control inputs to the model are acceleration ($\tilde{u}_1$), deceleration ($\tilde{u}_2$), and the power consumed for thermal management of the energy storage ($\tilde{u}_3$). The power to an external load ($\hat{u}$) should also be delivered by the energy storage. The power demand is assumed to be non-adjustable, e.g., an auxiliary load on the battery to fulfill a demand within the vehicle. 

The simplified vehicle model has the following structure:\vspace{-0.25cm}

\small
\begin{gather}\label{eq:vehModel}
\boldsymbol{x}_{k+1}=\begin{bmatrix}
1 & 1 & 0 & 0\\0 & 1 & 0 & 0\\0 & 0 & 1 & 0\\0 & 0 & 0 & 1
\end{bmatrix}\boldsymbol{x}_{k}~~~~~~~~~~~~~~~~~~~\\~~~~~~~~~~~~~~~+\begin{bmatrix}
1 & 1 & 0\\1 & -1 & 0\\-0.8 & 0.8 & -0.15\\1 & 1 & -0.85
\end{bmatrix}\tilde{\boldsymbol{u}}_{k}+\begin{bmatrix}
0\\0\\-0.25\\1
\end{bmatrix}\hat{{u}}_{k}\nonumber
\end{gather}
\normalsize

\hspace{-0.35cm}where $\boldsymbol{x}=[x_1,x_2,x_3,x_4]^\intercal$ and $\tilde{\boldsymbol{u}}=[\tilde{u}_1,\tilde{u}_2,\tilde{u}_3]^\intercal$. The controller design objective is to track a given desired position ($x_1^d$) while enforcing the state and input constraints. Specifically, the thermal index of the energy storage ($x_4$) is preferred to be maintained below its upper optimum operation limit $x_4^{UL}$. It is noted that while the vehicle tracks the reference position, the energy storage needs to deliver power to satisfy the external load ($\hat{u}$). Moreover, operating the vehicle by adjusting $\tilde{\boldsymbol{u}}$ and responding to $\hat{u}$ consumes the on-board vehicle energy. The operation of the vehicle is required to be optimized such that the on-board energy lasts until the end of the vehicle mission. The state and input constraints are listed below:\vspace{-0.15cm}
\begin{gather}\label{eq:vehConst_state}
[-1,-20,0,0]^\intercal\leq \boldsymbol{x}_{k}\leq [100,20,100,30]^\intercal\\
[-1,-1,-1]^\intercal\leq \tilde{\boldsymbol{u}}_k \leq [1,1,1]^\intercal. \label{eq:vehConst_input}
\end{gather}

In practice, the exact demand preview is assumed to be unknown over the entire vehicle mission, however, an approximate knowledge of the demand is available to the vehicle. Moreover, we assume the vehicle is connected to a server which can provide short-range and accurate demand predictions to the vehicle. The actual and approximate profiles of the power demand as the load ($\hat{u}$) on the on-board energy storage are shown in Fig.~\ref{fig:ErrorPredictionConcept}-(a,c). 

In order to evaluate the H-MPC performance for the vehicle model (\ref{eq:vehModel}), first a centralized single-layer MPC (S-MPC) is designed. To this end, we consider an S-MPC formulated over a finite-horizon ($N$) with $\tilde{u}_1,\tilde{u}_2,\tilde{u}_3$ being the optimization variables according to (\ref{eqn:S_MPC_gen}).~The stage cost of S-MPC, $\ell_{\text{S-MPC}}$, is defined as follows:\vspace{-0.2cm}
\begin{gather}\label{eqn:S_MPC_cost}
    \ell_{\textsuperscript{S-MPC}}=\lambda_1\tilde{u}_1^2+\lambda_2\tilde{u}_2^2+\lambda_3(x_1-x_1^d)^2
\end{gather}
where $x_1^d$ is the desired and known position trajectory. $\mathcal{X}$ and $\mathcal{U}$ are the convex sets properly defined according to~(\ref{eq:vehConst_state}) and (\ref{eq:vehConst_input}), respectively. The S-MPC optimization problem is solved at every time step, then the horizon is shifted by one step ($T=1~sec$), and only the current control is commanded to the system. The S-MPC simulation is carried out on a desktop computer, with an Intel\textsuperscript{\textregistered} Core i7@2.60 GHz processor, in MATLAB\textsuperscript{\textregistered}/SIMULINK\textsuperscript{\textregistered} using YALMIP~\cite{lofberg2004yalmip} for formulating the optimization problem, and IPOPT\cite{wachter2006implementation} for solving the optimization problem numerically. Additionally, in order to ensure the feasibility of the S-MPC solution, slack optimization variables are used to soften the constraints on the position ($x_1$) and thermal index ($x_4$). In all simulated cases, a very large weighting factor (e.g., $10^6$) has been considered for the slack variable term in the MPC stage cost to ensure a comparable performance. 

\begin{figure}[h!]
	\begin{center}
		\includegraphics[width=0.8\columnwidth]{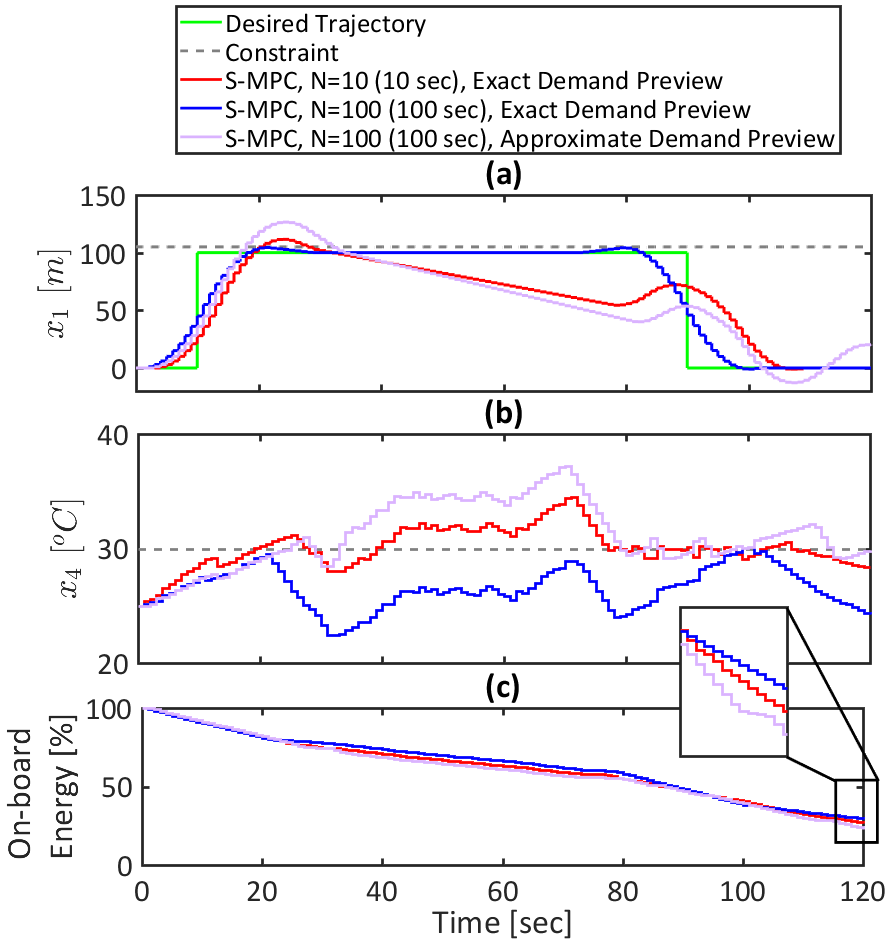} \vspace{-0.35cm}   		\caption{Performance of the S-MPC with different prediction horizon lengths: (a) tracking the desired position trajectory, (b) thermal index, and (c) on-board energy.}\vspace{-0.65cm} 
		\label{fig:HorizonPreviewImpact_SMPC_combine} 
	\end{center}
\end{figure}

Fig.~\ref{fig:HorizonPreviewImpact_SMPC_combine} shows the performing of the S-MPC with short ($N=10~(10~sec)$) and long ($N=100$) prediction horizons. Under the unrealistic assumption that the controller has the exact knowledge of the demand profile shown in Fig.~\ref{fig:ErrorPredictionConcept}-(c), Fig.~\ref{fig:HorizonPreviewImpact_SMPC_combine}-(a) shows that the S-MPC with long horizon tracks the desired position accurately with no $x_1$ constraint violation while enforcing the constraint on the thermal index state $x_4$. When the shorter horizon is used, the S-MPC cannot mitigate the violations of $x_4$ constraint; thus, it compromises $x_1$ tracking (Fig.~\ref{fig:HorizonPreviewImpact_SMPC_combine}-(b)) to avoid further increase in the thermal index. The reason for this deviation in the S-MPC performance with shorter prediction horizon is that the controller does not have enough lead time to take proactive actions in response to the external power demand early on. Moreover, compared to the long horizon S-MPC, the short horizon controller consumes more energy (Fig.~\ref{fig:HorizonPreviewImpact_SMPC_combine}-(c)).
\vspace{-0.3cm} 

It was shown that extending the prediction horizon improves S-MPC tracking and constraint satisfaction performance. While this observation is not a surprise, as discussed in Sec~\ref{sec:Sec2}, long horizon optimization for slow states is associated with demand forecast uncertainty and computational complexity. Fig.~\ref{fig:HorizonPreviewImpact_SMPC_combine} also shows the S-MPC performance with $N=100$ when the incorporated preview is based on the approximate knowledge of the demand profile (Fig.~\ref{fig:ErrorPredictionConcept}-(a)). As can be seen, the uncertainty in the demand preview significantly alters the long horizon S-MPC performance. The desired position is not being tracked preferably and the thermal index constraint violation occurs more often. Additionally, according to Fig.~\ref{fig:HorizonPreviewImpact_SMPC_combine}-(c), the S-MPC with the average demand preview has the highest energy consumption.

To address the demand forecast uncertainty impact on the S-MPC performance, a robust H-MPC is designed for the vehicle model (\ref{eq:vehModel}) according to the proposed framework in Secs.~\ref{sec:Sec2} and \ref{sec:sec3}. The scheduling layer MPC is based on the optimization problem formulated in~(\ref{eqn:H_MPC_UL_gen}) with the same stage cost defined in (\ref{eqn:S_MPC_cost}), and with a slower update rate.~The long prediction horizon of the scheduling layer MPC is set to $H_s=20~(100~sec)$. The sampling time of the scheduling layer MPC is $T_s=5~sec$, thus $\nu=5$. The robust constraint set $\mathcal{X}^{\text{\textbf{robust}}}$ is calculated according to the procedure (\ref{eq:X_robust_1})-(\ref{eq:X_robust_2}) proposed in Sec.~\ref{sec:sec3}. 

\begin{figure}[h!]
	\begin{center}
		\includegraphics[width=0.9\columnwidth]{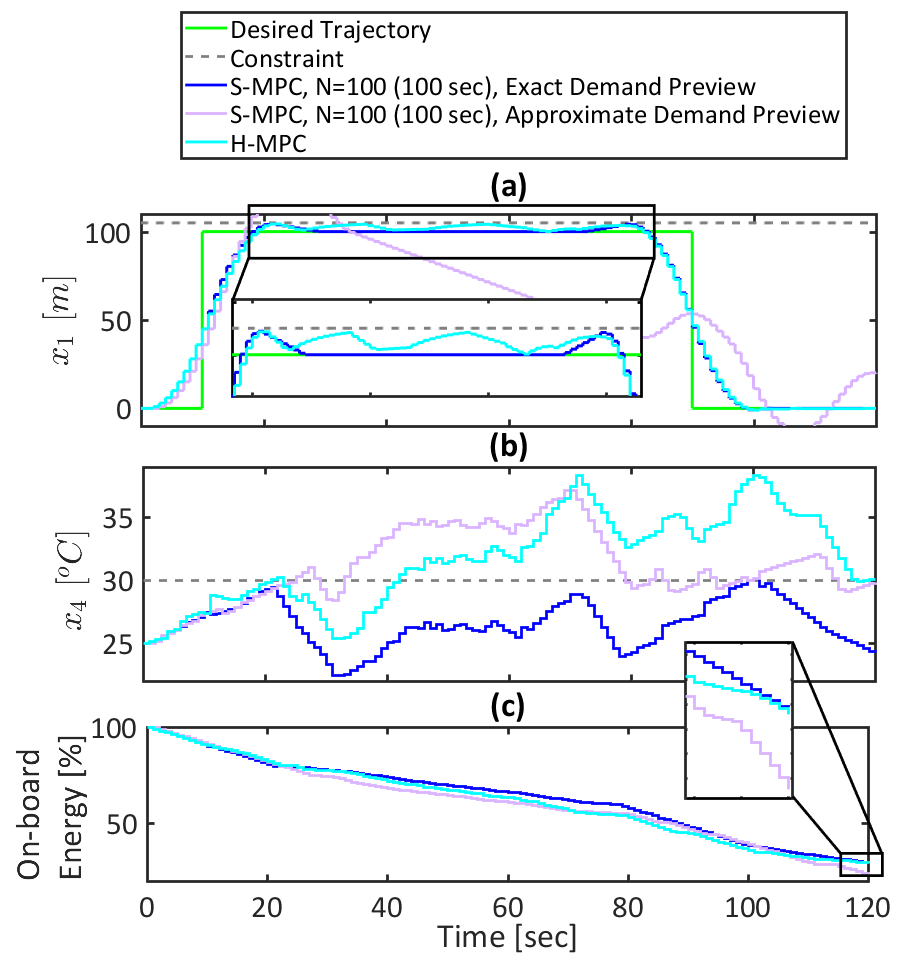} \vspace{-0.35cm}   
		\caption{Comparison between performances of S-MPC and H-MPC: (a) tracking the desired position trajectory, (b) thermal index, and (c) on-board energy ($H_p=20~sec$, $H_s=N=100~sec$).}\vspace{-0.7cm} 
		\label{fig:H_MPC_wTerminal_all} 
	\end{center}
\end{figure}

The scheduling layer MPC optimization problem is solved according to the average demand preview shown in~Fig.~\ref{fig:ErrorPredictionConcept}-(a). This solution is used to calculate the optimal thermal index trajectory ($\leftidx{^*}{x}_4^{\text{\textbf{slow}}}$). This scheduled $\leftidx{^*}{x}_4^{\text{\textbf{slow}}}$ trajectory is then incorporated in the stage cost of the piloting layer, which is formulated according to (\ref{eqn:H_MPC_LL_gen}) as a finite-horizon optimization problem over a relatively shorter prediction horizon ($H_p=20~(20~sec)$) and with a faster sampling period of $T_p=T=1~sec$. Note that it is assumed the exact knowledge of the demand preview is available to the piloting-layer MPC over the receding horizon $H_p=20~(20~sec)$.~The stage cost of the piloting layer MPC ($\ell_{\text{piloting}}$) is defined as:\vspace{-0.15cm}
\begin{gather} \label{eqn:H_MPC_LL_cost}
    \ell_{\text{piloting}}=\lambda_1\tilde{u}_1^2+\lambda_2\tilde{u}_2^2+\lambda_3(x_1-x_1^d)^2+\lambda_4(x_4-\leftidx{^*}{x}_4^{\text{\textbf{slow}}})^2.
\end{gather}
%
Since $T_s>T_p$, the scheduled $\leftidx{^*}{x}_4^{\text{\textbf{slow}}}$ trajectory is passed on as a piecewise constant function according to (\ref{eq:H_MPC_LL_Ref}). Fig.~\ref{fig:H_MPC_wTerminal_all} shows the comparison between the S-MPC and the baseline H-MPC with $\mathcal{X}$ as the constraint set. The H-MPC, unlike the S-MPC with average demand preview, manages to deliver an acceptable $x_1$ tracking performance (Fig.~\ref{fig:H_MPC_wTerminal_all}-(a)). Since the robust constraint set $\mathcal{X}^{\text{\textbf{robust}}}$ is not being incorporated in the baseline H-MPC (\ref{eqn:H_MPC_UL_gen}), often violation of thermal index state $x_4$ is still observed in Fig.~\ref{fig:H_MPC_wTerminal_all}-(b). Additionally, compared to the S-MPC with long prediction horizon and average demand preview, the H-MPC shows less energy consumption (Fig.~\ref{fig:H_MPC_wTerminal_all}-(c))~
Overall, the H-MPC slightly relaxes the accurate $x_1$ tracking requirement by putting part of the effort into tracking the scheduled $\leftidx{^*}{x}_4^{\text{\textbf{slow}}}$. However, due to the uncertainty in the long-range demand forecast used at the scheduling layer, $x_4$ constraint enforcement is not achieved. Note that operating the vehicles while the thermal index is maintained close to the upper limit $x_4^{UL}$ is equivalent to less energy consumption. Since the demand preview is uncertain, steering $x_4$ to $x_4^{UL}$ (e.g., $30^oC$) increases the chances of constraint violation once the vehicle faces the actual power demand. 

\begin{figure}[h!]
	\begin{center}
		\includegraphics[width=0.9\columnwidth]{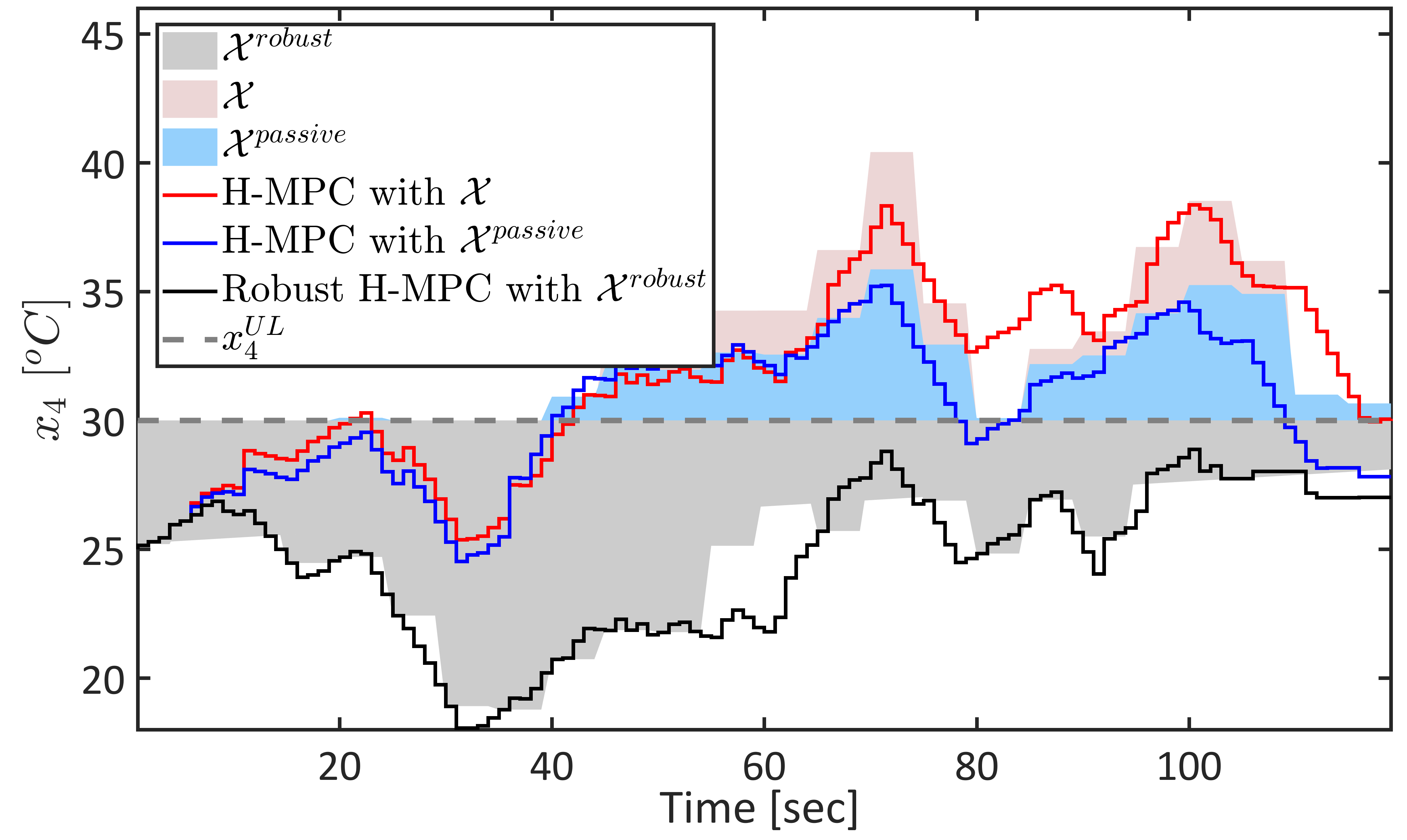} \vspace{-0.45cm} \caption{Comparison between performances of H-MPC and robust H-MPC in regulating the energy storage thermal index .}\vspace{-0.5cm} 
		\label{fig:R_H_MPC_wTerminal_Tbat} 
	\end{center}
\end{figure}

Upon incorporation of $\mathcal{X}^{\text{\textbf{robust}}}$ in (\ref{eqn:H_MPC_UL_gen}), the robust H-MPC (\ref{eqn:H_MPC_UL_gen_robust}) becomes able to take advantage of the short horizon and more accurate demand preview over the piloting layer horizon. As shown in Fig.~\ref{fig:R_H_MPC_wTerminal_Tbat}, compared to the baseline H-MPC (\ref{eqn:H_MPC_UL_gen}) with $\mathcal{X}$ as the constraint set at the scheduling layer, the robust H-MPC (\ref{eqn:H_MPC_UL_gen_robust}) with $\mathcal{X}^{\text{\textbf{robust}}}$ shows significantly better $x_4$ constraint enforcement under the uncertainty associated with long horizon demand forecast. As shown in Fig.~\ref{fig:R_H_MPC_wTerminal_Tbat}, the robust H-MPC starts to tighten the $x_4^{UL}$ constraint early on to ensure its enforcement during the period that the demanded load on the battery is higher and any uncertainty in the demand preview could lead to violation of the thermal index limit.\vspace{-0.15cm} 

The performance of the robust H-MPC is also compared with the H-MPC with $\mathcal{X}^{\text{\textbf{passive}}}$, which is computed according to (\ref{eq:X_passive}). Unlike $\mathcal{X}^{\text{\textbf{robust}}}$, in the passive approach $\mathcal{X}$ is tightened whenever the violation of $x_4^{UL}$ is detected by the piloting layer MPC. This means the scheduling layer MPC is informed about the constraint violation only after its occurrence, then it manages to update the scheduled trajectories to reduce the chances of constraint violation in subsequent time steps. Fig.~\ref{fig:R_H_MPC_wTerminal_Tbat} shows that compared to the baseline H-MPC, the H-MPC with $\mathcal{X}^{\text{\textbf{passive}}}$ can only slightly decrease the $x_4$ constraint violation. This can be explained with respect to the slower dynamics of $x_4$ which call for a proactive constraint tightening algorithm to effectively improve the H-MPC robustness. This requirement has been addressed by the proposed robust H-MPC framework in this paper. Note that with baseline and passive H-MPCs, enforcement of the original $x_4^{UL}$ constraint is infeasible. As a results, the computed $\mathcal{X}$ and $\mathcal{X}^{\text{\textbf{passive}}}$ shown in Fig.~\ref{fig:R_H_MPC_wTerminal_Tbat} are based on the softened $x_4$ constraint using slack variables.

\begin{figure}[t!]
	\begin{center}
		\includegraphics[width=0.8\columnwidth]{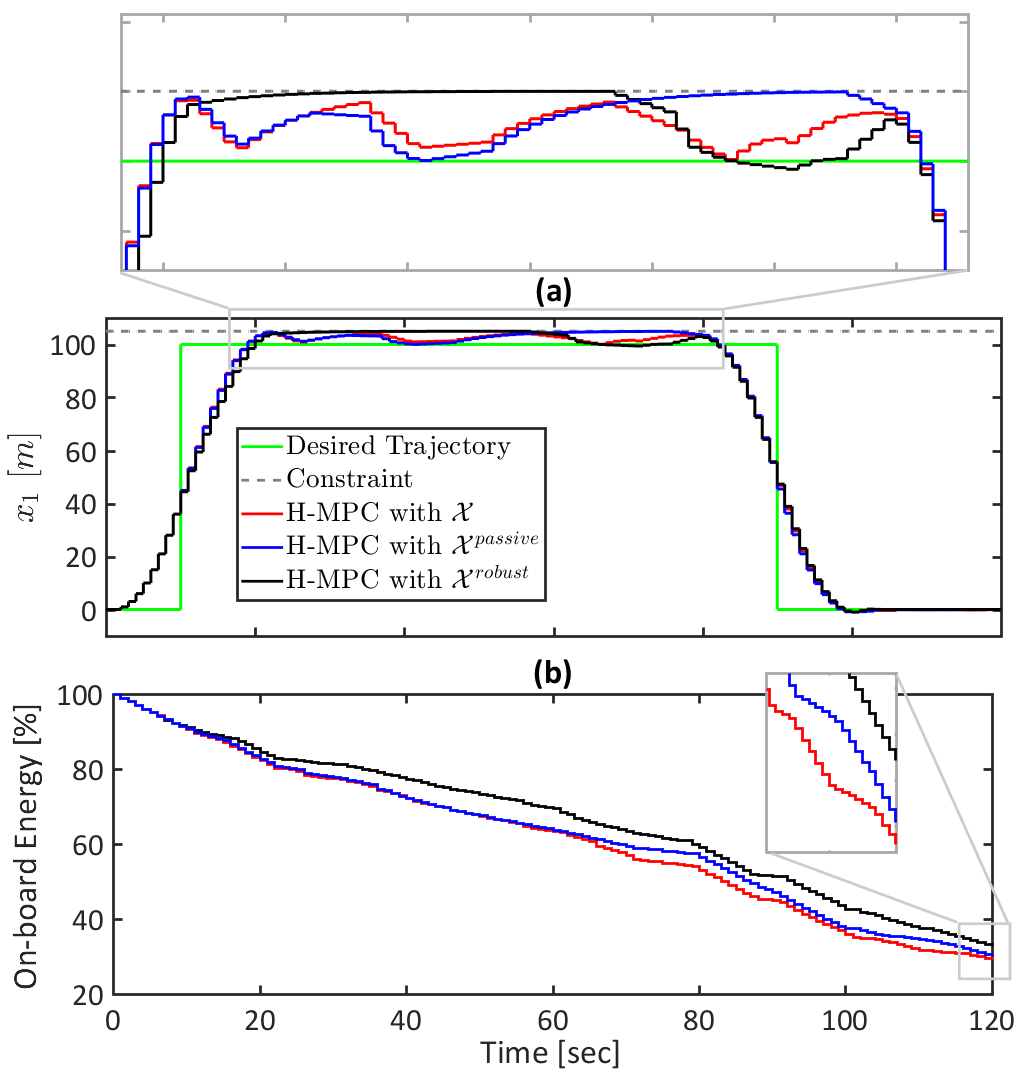} \vspace{-0.45cm}   
		\caption{Comparison between performances of H-MPC and robust H-MPC: (a) tracking the desired position trajectory, and (b) on-board stored energy ($H_p=20~sec$, $H_s=N=100~sec$)}\vspace{-0.95cm} 
		\label{fig:R_H_MPC_wTerminal_all} 
	\end{center}
\end{figure}

The other interesting observation about the robust H-MPC is its performance in tracking $x_1^d$ (Fig.~\ref{fig:R_H_MPC_wTerminal_all}-(a)) while enforcing $x_4^{UL}$ constraint. Given the limited energy stored on-board of the vehicle, the robust H-MPC slightly relaxes the $x_1$ trajectory tracking objective (while it enforces $x_1$ constraint) during the period it puts more effort to mitigate $x_4^{UL}$ violation by decreasing (cooling) the energy storage thermal index from around $t=10~sec$ to $t=60~sec$. After this time period, since the robust H-MPC does not see any other major external power demand over the long prediction horizon, it increases the position tracking priority, see Fig.~\ref{fig:R_H_MPC_wTerminal_all}.

The energy consumption results of different studied H-MPCs are shown in Fig.~\ref{fig:R_H_MPC_wTerminal_all}-(b). It can be seen that the robust H-MPC has the lowest energy consumption, as compared to the baseline and passive H-MPCs. This is an interesting observation as it shows the robust H-MPC framework does not necessarily lead to a conservative controller design which is usually expected when designing a robust MPC under the influence of unknown external disturbances. We recall that for the specific class of dynamic systems with slow states considered in this paper, steering the slow state to its limit usually leads to lower energy consumption. The constraint violation often occurs when the external power demand on the energy storage is higher, meaning during the time periods that the energy storage is required to deliver power for operation of the vehicle, it also has to provide extra power for regulating the thermal state. The robust H-MPC, on the other hand, effectively shifts the thermal management load to those periods during which the power demand is relatively lower. Thanks to this intelligent thermal management load shift by the robust H-MPC, not only less often thermal constraint violations are observed, but the overall efficiency (fuel economy) of the vehicle is also higher, compared to the baseline H-MPC. 

\vspace{-0.1cm}
\section{Summary and Conclusions}\label{sec:5}
The problem of robust hierarchical MPC (H-MPC) design for constrained systems with multi-timescales and exposed to external demands was considered in this paper. From a new prospective, we proposed to use demand preview information with different accuracies at different layers of an H-MPC to proactively predict the (chances of) constraint violations and compute a robust constrain set to enhance the overall robustness of the H-MPC against the uncertainty associated with long horizon demand forecast. The computation of the robust constraint set is built upon the information that is available at each layer of an H-MPC, but they are not often being communicated between the layers. By leveraging this information and accessing to demand preview information with different accuracies, we showed that the robustness of conventional H-MPC can be effectively improved. We demonstrated the application of the proposed robust H-MPC for a simplified vehicle thermal management case study. The simulation results confirmed the effectiveness of the proposed approach in regulating the thermal state of the vehicle based on an uncertain knowledge of the demand preview while improving the overall system efficiency, as compared to a conventional H-MPC framework.

\vspace{-0.15cm}

\bibliographystyle{unsrt} 
\bibliography{ACC2018Ref.bib} 

\end{document}